\documentclass[final,12pt]{colt2026} 

\usepackage{times}
\usepackage{algorithm,algorithmic}

\usepackage{amsmath,amsfonts,bm,amssymb}

\def\eqref#1{equation~\ref{#1}}

\def\1{\bm{1}}

\def\vzero{{\bm{0}}}

\def\va{{\bm{a}}}
\def\vb{{\bm{b}}}
\def\vc{{\bm{c}}}

\def\vg{{\bm{g}}}
\def\vh{{\bm{h}}}

\def\vu{{\bm{u}}}
\def\vv{{\bm{v}}}

\def\vx{{\bm{x}}}
\def\vy{{\bm{y}}}
\def\vz{{\bm{z}}}

\def\mA{{\bm{A}}}

\def\mH{{\bm{H}}}
\def\mI{{\bm{I}}}

\def\mO{{\bm{O}}}

\def\mQ{{\bm{Q}}}

\def\mLambda{{\bm{\Lambda}}}

\def\mLambda{{\bm{\Lambda}}}

\DeclareMathAlphabet{\mathsfit}{\encodingdefault}{\sfdefault}{m}{sl}
\SetMathAlphabet{\mathsfit}{bold}{\encodingdefault}{\sfdefault}{bx}{n}

\def\gO{{\mathcal{O}}}

\def\gX{{\mathcal{X}}}

\def\sB{{\mathbb{B}}}

\def\sN{{\mathbb{N}}}
\def\sO{{\mathbb{O}}}

\def\sR{{\mathbb{R}}}

\DeclareMathOperator*{\argmin}{arg\,min}

\newtheorem{thm}{Theorem}[section]
\newtheorem{dfn}{Definition}[section]

\newtheorem{lem}{Lemma}[section]
\newtheorem{asm}{Assumption}[section]

\newtheorem{prop}{Proposition}[section]

\usepackage{cancel}
\usepackage{nicefrac}
\usepackage{makecell}
\usepackage{tcolorbox}
\usepackage{booktabs}%
\usepackage{bbding}
\usepackage{pifont}
\definecolor{mydarkgreen}{RGB}{39,130,67}
\definecolor{mydarkred}{RGB}{192,25,25}

\title[Faster Newton Methods for Convex and Nonconvex Optimization in Gradient Complexity]{Faster Newton Methods for Convex and Nonconvex Optimization in Gradient Complexity}

\coltauthor{%
 \Name{Lesi Chen} \Email{chenlc23@mails.tsinghua.edu.cn}\\
 \addr IIIS, Tsinghua University\\
       Shanghai Qizhi Institute
 \AND
 \Name{Chengchang Liu} \Email{ccliu22@cse.cuhk.edu.hk} \\
 \addr 
 Department of Computer Science $\&$ Engineering, The Chinese University of Hong Kong 
 \AND 
 \Name{Luo Luo} \Email{luoluo@fudan.edu.cn} \\
 \addr {School of Data Science, Fudan University  \\
 Shanghai Key Laboratory for Contemporary Applied Mathematics}
 \AND 
   \Name{Jingzhao Zhang}${}^\dagger$ 
 \Email{jingzhaoz@mail.tsinghua.edu.cn} \\
 \addr {IIIS, Tsinghua University\\
        Shanghai AI Lab\\
       Shanghai Qizhi Institute
 }
}

\begin{document}

\maketitle

\begin{abstract}%
Second-order optimization methods are computationally expensive for large-scale problems. 
Recently, Doikov, Chayti, and Jaggi (ICML 2023) proposed the LazyCRN method that reduces computation by studying the gradient complexity of second-order methods. Their method can achieve a gradient complexity of $\mathcal{O}( \bar d + \bar d^{1/2} \epsilon^{-3/2})$ and  $\mathcal{O}( \bar d + \bar d^{1/2} \epsilon^{-1/2})$ for nonconvex and convex optimization, respectively, where $\bar d$ is the effective dimension and $\epsilon$ is the target precision. 
Very recently, Adil, Bullins, Sidford, and Zhang (NeurIPS 2025) improved the gradient complexity to $\mathcal{O}( \bar d + \bar d^{1/3} \epsilon^{-3/2} \ln^{18} \epsilon^{-1})$ for nonconvex optimization. However, the tightness of these methods remains open.
In this work, we propose new methods that
achieve an improved complexity of $\mathcal{O}( \bar d + \bar d^{1/3} \epsilon^{-3/2})$ and  $\mathcal{O}( (\bar d + \bar d^{13/21} \epsilon^{-2/7}) \ln \bar d)$ for nonconvex and convex optimization, respectively, improving best-known results for both setups.
\end{abstract}

\section{Introduction}

Minimizing differentiable functions
has wide applications in machine learning, including matrix completion \citep{hardt2014understanding}, phase retrieval \citep{candes2015phase}, and neural network training \citep{lecun2015deep}. 
In this work, we focus on second-order methods \citep[Chapter 4]{nesterov2018lectures}, which achieve faster convergence rates compared to their first-order counterparts. Over the past two decades,
optimal second-order methods measured by oracle complexities \citep{nemirovskij1983problem} (\textit{i.e.}, the number of gradient and Hessian calls) have been proposed. Remarkably, for nonconvex problems, \citet{nesterov2006cubic} proposed the optimal second-order method known as cubic regularized Newton (CRN). For convex problems, \citet{monteiro2013accelerated} proposed a near-optimal method known as accelerated Newton proximal extragradient (A-NPE), and recently
\citet{carmon2022optimal,kovalev2022first} proposed optimal A-NPE methods.

However, oracle-optimal methods could be slow in practice due to the computational burden of solving auxiliary problems. An example is cutting-plane methods: while they have long been known to have optimal first-order oracle complexity \citep{nemirovskij1983problem}, they often require solving non-trivial subproblems, and thus a wealth of research has been dedicated to making them more computationally efficient \citep{vaidya1989new,vaidya1996new,lee2015faster,jiang2020improved} over the past three decades.
A similar situation holds for second-order algorithms: optimal second-order methods could also be costly due to the computation of the Hessian oracle and matrix inversion. 
It motivates us to better model computation for designing faster second-order methods.

Recently, \citet{doikov2023second} proposed to measure the complexity of second-order methods by 
\textit{equivalent gradient calls}. This new computation model takes into account that Hessian oracles are usually more expensive than gradient oracles in practice, and 
thus treats a single Hessian query as equivalent to $\bar d$ gradient queries, where $\bar d$ is the \textit{effective dimension}. Additionally, for second-order methods, the equivalent gradient complexity can be easily translated into actual computational complexity (see Section \ref{subsec:grad-complexity}). 
Subsequently, new second-order methods with low \textit{equivalent gradient complexities} are proposed \citep{doikov2023second,chayti2023unified,adil2025balancing}.

For nonconvex problems, 
\citet{doikov2023second} proposed the lazy cubic regularized Newton (LazyCRN) method
and showed that it can find an $\epsilon$-stationary point a nonconvex function in an equivalent gradient complexity of $\gO(\bar d + \bar d^{1/2} \epsilon^{-3/2})$, which improves over the 
$\gO(\bar d \epsilon^{-3/2})$ complexity by original CRN method \citep{nesterov2006cubic}. 
Very recently, \citet{adil2025balancing} proposed a novel restarted approximate Hessian accelerated gradient descent (RAH-AGD)  method, achieving an improved complexity of 
$\gO(\bar d + \bar d^{1/3} \epsilon^{-3/2} \ln^{18} \epsilon^{-1})$ by applying the restarted AGD \citep{li2023restarted} preconditioned with a carefully designed approximate lazy Hessian. However, the $\epsilon$-dependency in the complexity bound of RAH-AGD involves an extra $\ln^{18} \epsilon^{-1}$ factor compared to the $\Omega(\epsilon^{-3/2})$
lower bound \citep{carmon2020lower}, suggesting room for improvement. 


For convex problems, \citet{chayti2023unified} showed that LazyCRN can find an $\epsilon$-global solution in the complexity of  $\gO(\bar d + \bar d^{1/2} \epsilon^{-1/2})$, which also improves the 
$\gO(\bar d \epsilon^{-1/2})$ complexity by CRN \citep{nesterov2006cubic}. However, the 
dependency in $\epsilon$ is far from the optimal $\Omega(\epsilon^{-2/7})$ lower bound~\citep{arjevani2019oracle}.
The acceleration of LazyCRN for convex problems has been regarded as an open problem in 
in the prior literature \citep{doikov2023second,doikov2025first}.

\paragraph{Our Contributions.} Inspired by known results, we improve gradient complexities for both convex and nonconvex problems by proposing new second-order methods:
\begin{enumerate}
    \item For nonconvex problems, we propose a novel second-order method called Nonconvex Accelerated Lazy Extra Newton (NALEN).
    We show that NALEN can provably find an $\epsilon$-stationary point in an equivalent gradient complexity of $\gO( \bar d+ \bar d^{1/3} \epsilon^{-3/2} )$. Compared with RAH-AGD \citep{adil2025balancing}, our algorithm is easy to analyze and removes the $\ln^{18} \epsilon^{-1}$ factor in the complexity bound.
    \item For convex problems, we also propose an efficient second-order method called Convex Accelerated Lazy Extra Newton (CALEN) by using a restarted NALEN algorithm to solve the proximal subproblems in the optimal second-order methods \citep{carmon2022optimal}.
    We show that CALEN can find an $\epsilon$-global solution in an equivalent gradient complexity of 
    $ \gO \left( (\bar d + \bar d^{13/21} \epsilon^{-2/7})\ln \bar d
    \right)$, which improves the prior art of $\gO(\bar d+ \bar d^{1/2} \epsilon^{-1/2})$ achieved by LazyCRN \citep{chayti2023unified}. 
\end{enumerate}
Our results significantly improve the known ones for both nonconvex and convex optimization, as compared in Table \ref{tab:main-result-nonconvex} and \ref{tab:main-result-convex}, respectively.

\begin{table*}[htbp] 
\caption{We compare the equivalent gradient complexities of different methods to find an $\epsilon$-stationary point of $d$-dimensional nonconvex problems (\ref{prob:main}) when the objective has Lipschitz continuous Hessians, where $\bar d \le d$ is the effective dimension.
}
\label{tab:main-result-nonconvex}
\centering
\begin{tabular}{c c c c }
\hline 
Order & Method  & Reference & Equivalent Gradient Complexity \\
\hline \hline  \addlinespace
& GD & \citet{carmon2020lower} & $\gO \left( \epsilon^{-2} \right)$  \\ \addlinespace 
1st & R-AGD & \citet{li2023restarted} &  $\gO \left( \epsilon^{-7/4} \right)$  \\ \addlinespace
& OQN & \citet{jiang2025improved} & $ \gO\left(  d^{1/4}  
\epsilon^{-13/8} \right)$ \\ \addlinespace
& Lower Bound & \citet{carmon2021lower} & $\Omega(\epsilon^{-12/7})$ \\ \addlinespace
\hline \addlinespace
& CRN &  \citet{nesterov2006cubic} &  $ \gO \left( \bar d \epsilon^{-3/2}  \right)$  \\ \addlinespace
& LazyCRN & \citet{doikov2023second} & $ \gO \left( \bar d + \bar d^{1/2} \epsilon^{-3/2} \right)$ \\
\addlinespace
2nd & RAH-AGD & \citet{adil2025balancing} & 
$\gO \left( \bar d + \bar d^{1/3} \epsilon^{-3/2} \ln^{18} \epsilon^{-1}  \right)$ 
\\ \addlinespace
& ~~~~NALEN (ours)  & Theorem \ref{thm:ALEN-NC} 
& $ \gO \left( \bar d + \bar d^{1/3} \epsilon^{-3/2}
\right)$  
\\ \addlinespace
& Lower Bound & \citet{carmon2020lower} & $\Omega(\epsilon^{-3/2})$ \\ \addlinespace
\hline 
\end{tabular}
\end{table*}

\begin{table*}[htbp] 
\caption{We compare the equivalent gradient complexities of different methods to find an $\epsilon$-solution of  $d$-dimensional convex problems (\ref{prob:main}) when the objective has Lipschitz continuous Hessians, where $\bar d \le d$ is the effective dimension.
}
\label{tab:main-result-convex}
\centering
\begin{tabular}{c c c c }
\hline 
Order & Method  & Reference & Equivalent Gradient Complexity  \\
\hline \hline  \addlinespace
& FastCPM  &  \citet{lee2015faster} & $ \gO \left( d \ln \epsilon^{-1} \right)$ \\ \addlinespace
1st & AGD & \citet{nesterov1983method} & $\gO \left( \epsilon^{-1/2} \right)$  \\ \addlinespace 
& A-QNPE &
\citet{jiang2024accelerated}
& $\gO \left( 
d^{1/5}  \epsilon^{-2/5} \ln (d\epsilon^{-2})
\right)$ 
 \\ \addlinespace
& Lower Bound & \citet{nemirovskij1983problem} & $\Omega(\epsilon^{-1/2})$ \\
\hline \addlinespace
& OptANPE & \citet{carmon2022optimal} & $ \gO \left( \bar d \epsilon^{-2/7} \right)$ \\
\addlinespace
2nd & LazyCRN & \citet{chayti2023unified} & 
$\gO \left( \bar d + \bar d^{1/2} \epsilon^{-1/2}  \right)$ 
\\ \addlinespace
& ~~~CALEN (ours) & Theorem \ref{thm:ALEN}
& $ \gO \left( (\bar d + \bar d^{13/21} \epsilon^{-2/7}) \ln \bar d
\right)$  
\\ \addlinespace
& Lower Bound & \citet{arjevani2019oracle} & $\Omega(\epsilon^{-2/7})$ \\
\hline 
\end{tabular}
\end{table*}

\paragraph{Notations.} We use $\Vert \cdot \Vert$ to denote the Euclidean norm of a vector and the spectral norm of a matrix. $\sB_D^d(\vc)$ denotes the ball $ \{  \vx \in \sR^d : \Vert \vx - \vc \Vert \le D \}$. For any set $S$ and functions $g, h :  S \rightarrow [0, \infty)$ we write $g = \gO(h)$ or $h = \Omega(g)$ equivalently if
there exists $c > 0$ such that $g(s) \le c  h(s)$ for every $s \in S$. We also write 
$h = \Theta(g)$ if $h = \gO(g)$ and $h = \Omega(g)$ hold simultaneously.

\section{Additional Related Works}


\subsection{Convex optimization}

For first-order methods,
\citet{nesterov1983method} proposed the accelerated gradient descent (AGD) that can find an $\epsilon$-solution to a convex problem with $\gO(\epsilon^{-1/2})$ gradient calls, which is oracle-optimal as it matches the $\Omega(\epsilon^{-1/2})$ lower bound \citep{nemirovskij1983problem}. 
In low dimensions, the cutting plane methods (CPM) typically enjoy better theoretical upper bounds and also achieve the optimal gradient complexity of $\gO(d \ln \epsilon^{-1})$ \citep{nemirovskij1983problem,braverman2020gradient}. A variety of works such as
\citep{vaidya1989new,vaidya1996new,lee2015faster,jiang2020improved} further improved the computational costs in auxiliary problems in CPMs. 

For second-order methods,
\citet{nesterov2006cubic} proposed the cubic regularized Newton (CRN) method that requires a $\gO( \bar d \epsilon^{-1/2})$ equivalent gradient complexity.
Furthermore, \citet{nesterov2008accelerating} proposed the accelerated cubic regularized Newton (A-CRN) method that achieves a better complexity of $\gO( \bar d \epsilon^{-1/3})$. Later,
\citet{monteiro2013accelerated} proposed the accelerated Newton proximal extragradient (A-NPE) method that further improves the complexity to $\gO(\bar d \epsilon^{-2/7} \ln \epsilon^{-1})$,
which is refined to 
$\gO(\bar d \epsilon^{-2/7})$ independently by \citet{carmon2022optimal,kovalev2022first}, whose $\epsilon$ dependencies are optimal as they match the lower bound provided by
\citet{arjevani2019oracle}. 
Additionally, 
\citet{jiang2024accelerated} proposed the accelerated quasi-Newton proximal extragradient (A-QNPE) method with a gradient complexity of $ \gO \left( 
d^{1/5} \epsilon^{-2/5}  \ln (d \epsilon^{-2})
\right)$.

\subsection{Nonconvex Optimization}
For first-order nonconvex optimization, it is well-known that gradient descent (GD) finds an $\epsilon$-stationary point for smooth functions in an $\gO(\epsilon^{-2})$ gradient complexity, which matches the lower bounds 
in \citep{carmon2020lower,carmon2021lower,chewi2023complexity}.
Faster iteration complexities can be achieved by exploiting higher-order smoothness and even using higher-order oracles:

For second-order methods, 
\citet{nesterov2006cubic} proposed the CRN algorithm that achieves an
$\gO( \bar d \epsilon^{-3/2})$ gradient complexity for second-order smooth functions
and subsequently \citet{carmon2020lower} proved matching lower bounds on the $\epsilon$ dependency.
In a novel work,
\citet{doikov2023second} proposed the lazy cubic regularized Newton (LazyCRN) method that improves the complexity of CRN to $\gO \left( \bar d + \bar d^{1/2} \epsilon^{-3/2}  \right)$, and subsequently
\citet{chayti2023unified} provided extensions to convex and gradient-dominated functions.
Very recently,
\citet{jiang2025improved} proposed an optimistic quasi-Newton (OQN) method with a complexity of $\gO( d^{1/4}  
\epsilon^{-13/8} 
)$.
\citet{adil2025balancing} proposed the RAH-AGD method that achieves a complexity of $\gO \left( \bar d +\bar d^{1/3} \epsilon^{-3/2} \ln^{18} \epsilon^{-1} \right)$.  

Under additional second-order smoothness assumptions,
\citet{carmon2018accelerated,agarwal2017finding} 
proposed novel accelerated first-order methods that achieve a $ \gO( \epsilon^{-7/4}\ln (d\epsilon^{-1})  + \bar d \epsilon^{-3/2})$ complexity by solving the CRN subproblem with AGD, which is further refined to 
$\gO( \epsilon^{-7/4} \ln \epsilon^{-1})$ by \citet{carmon2017convex,jin2018accelerated}, and more recently to $\gO( \epsilon^{-7/4} )$ by \citet{li2023restarted}.

For third-order smooth problems, \citet{carmon2017convex} also established an $\gO( \epsilon^{-5/3} \ln \epsilon^{-1})$ complexity  using the nonconvex AGD algorithm, while 
\citet{hinder2018cutting} showed a $ \gO( d^2 \epsilon^{-4/3} \ln \epsilon^{-1} )$ complexity by extending CPM to nonconvex optimization.
\citet{carmon2021lower} established nearly tight lower bounds for first-order methods under high-order smoothness assumptions.

\section{Preliminaries} \label{sec:pre}

This section introduces the problem setup and complexity measurements in this paper.

\subsection{Problem Setup}

This paper considers the optimization problem
\begin{align} \label{prob:main}
\min_{\vx \in \sR^d} f(\vx).
\end{align}
As a standard assumption for second-order methods \citep[Chapter 4]{nesterov2018lectures}, we assume that the objective function in Problem (\ref{prob:main}) has Lipschitz continuous Hessians.

\begin{asm} \label{asm:Hessian-Lip}
    We assume that $\nabla^2 f(\vx)$ exists and is $L$-Lipschitz continuous ($L>0$), \textit{i.e.},
    \begin{align*}
        \Vert \nabla^2 f(\vx) - \nabla^2 f(\vx') \Vert \le L \Vert \vx-  \vx' \Vert, \quad \forall \vx,\vx' \in \sR^d.
    \end{align*}
\end{asm}
We consider two typical problem classes characterized by the following two assumptions on the function $f:\sR^d \rightarrow \sR$, the first one is convex problems with nonempty solution sets and the second is nonconvex problems with lower-bounded objective values.



\begin{asm} \label{asm:convex}
We assume that $f$ is convex and has at least one solution $\vx^* \in \arg \min_{\vx \in \sR^d} f(\vx)$.
\end{asm}

\begin{asm} \label{asm:nonconvex}
We assume that $f$ can be nonconvex and we have $f^* := \inf_{\vx \in \sR^d} f(\vx) > - \infty$.
\end{asm}



Next, we introduce the optimality criteria for these problem classes. For convex problems, we target to find an approximate global solution that has a small suboptimality gap \citep{nesterov2018lectures}. For nonconvex problems, finding global solutions may require exponential time, and we target to find an approximate stationary point as a common compromise \citep{lan2020first,carmon2020lower}.

\begin{dfn} \label{dfn:sub-opt-point}
We call $\hat \vx \in \sR^d$  an $\epsilon$-(global) solution to Problem (\ref{prob:main}) if $f(\hat \vx) - f^* \le \epsilon$.
\end{dfn}

\begin{dfn} \label{dfn:sta-point}
We call $\hat \vx \in \sR^d$ an $\epsilon$-stationary point to Problem (\ref{prob:main}) if $\Vert \nabla f(\hat \vx) \Vert \le \epsilon$.
\end{dfn}

Finally, we need to define our complexity measurement. In the literature, it is common to measure the complexity of an algorithm by the number of oracle calls \citep{nemirovskij1983problem,nesterov2018lectures}. However, this framework does not distinguish the different costs of Hessian and gradient evaluations. Consequently, we then resort to the following notion of equivalent gradient complexity.

\subsection{Equivalent Gradient Complexity} \label{subsec:grad-complexity}

In this work, we follow \citep{doikov2023second,adil2025balancing} to unify the complexity measurements of first- and second-order methods by the number of \textit{equivalent gradient} calls.

\begin{dfn}
Let $c_\mH$ and $c_\vg$ be the computational costs of Hessian and gradient evaluations, respectively.
We define the effective dimension $\bar d$ as $\bar d := c_{\mH} / c_{\vg}$.
\end{dfn}

As mentioned in \citep{doikov2023second}, we assume that $\bar d \le d$ without loss of generality, since we can always compute a Hessian by the Hessian-vector-products in $d$ basics, and a Hessian-vector product can be further approximated to any precision by the difference between two gradients. Equipped with this definition, the costs of oracle evaluations in an algorithm are
\begin{align} \label{eq:total-oracle-cost}
    \text{Total oracle costs} = c_\mH n_\mH + c_\vg  n_\vg = c_\vg \left( \bar d n_\mH + n_\vg \right),
\end{align}
where $n_\mH$ and $n_\vg$ are the number of Hessian and gradient calls, respectively.
To minimize the total oracle costs in Eq. (\ref{eq:total-oracle-cost}), it is clear that for any $\bar d \ge 1$ the optimal setting is $n_\mH = n_\vg / \bar d$. \citet{doikov2023second} propose a natural way to achieve this optimal trade-off, that is, querying a new Hessian every~$m$ iterations and reusing the Hessian of the snapshot point. 

According to Eq. (\ref{eq:total-oracle-cost}), the optimal setting of $m$ to minimize the total oracle costs is $m = \Theta(\bar d)$.
In Section \ref{sec:imple}, we will also show that, if further considering the computational costs of solving Newton subproblems in second-order methods, the optimal choice of $m$ is $m = \Theta( (c_\mH + d^\omega ) / (c_\vg+ d^2) ) $, where $d^\omega \approx d^{2.37}$ is the asymptotic time for multiplying two $d \times d$ matrices \citep{duan2023faster}. For the convenience of describing this \textit{lazy Hessian} update, we follow \citet{doikov2023second} to introduce the following additional notations.

\paragraph{Additional Notations.}
Throughout this paper, we let $m$ be the period of query new Hessians (whose optimal value is fixed), and let $\pi(t) = t - (t \mod m)$ be the snapshot.    

\section{Results for Nonconvex Optimization}

We propose the NALEN method for nonconvex optimization. This method consists of two loops, where the outer loop applies the online-to-nonconvex conversion (O2NC) \citep{cutkosky2023optimal} to reduce the nonconvex optimization problem to learning a descent direction, and the inner loop applies a modified Lazy Extra-Newton (LEN) algorithm \citep{chen2025secondorder} to solve the subproblem.

\subsection{Outer Loop: Online-to-Nonconvex Conversion} \label{subsec:back-o2nc}

\citet{cutkosky2023optimal} reformulated the task of finding stationary points in nonconvex optimization into the following online learning problem of finding a descent direction.

\begin{tcolorbox}[colback=white, colframe=black, boxrule=0.2mm] 
\textbf{Online Direction Learning Problem} \vspace{1mm} \\ 
For $n = 0 ,\cdots, KT-1$
\begin{enumerate}
    \item The learner chooses $\Delta_{n+1}  \in \sB_D^d(\vzero)$. 
    \item The environment computes
    \begin{align} \label{eq:update-w-x}
     \vy_{n+1} = \vx_n + \nicefrac{\Delta_{n+1}}{2} \quad {\rm and} \quad \vx_{n+1} = \vx_n + \Delta_{n+1}.
    \end{align}
    \item The learner observes loss $\langle \nabla f(\vy_{n+1}), \Delta_{n+1} \rangle $.
\end{enumerate}
\textbf{Goal:} Minimize the the $K$-shifting regret:
\begin{align} \label{dfn:shifting-reg}
    R_T\left(\vu^{(1)},\cdots,\vu^{(K)} \right):= \sum_{k=1}^K \sum_{n = (k-1)T +1}^{kT} \langle \nabla f (\vy_n) , \Delta_n - \vu^{(k)}\rangle,
\end{align}
where $\vu^{(k)}$ is the average descent direction in the $k$th epoch defined as 
\begin{align} \label{dfn:uk}
    \vu^{(k)} = - \nicefrac{D \sum_{n = (k-1) T +1}^{kT} \nabla f(\vy_n)}{\Vert \sum_{n = (k-1) T +1}^{kT} \nabla f(\vy_n) \Vert}
\end{align}
\end{tcolorbox}

Intuitively, a learner who can minimize this regret must have sufficient knowledge of the descent directions $\vu^{(k)}$, which suffices to identify an approximate stationary point. This intuition can be formalized in the following lemma, which is first implicitly used by \citet{cutkosky2023optimal}, and recently explicitly stated by \citet{jiang2025improved}.

\begin{lem}[{\citet[Proposition 2.3]{jiang2025improved}}] \label{lem:o2nc}
Under Assumption \ref{asm:Hessian-Lip}, we have
\begin{align*} 
\frac{1}{K} \sum_{k=1}^K \left \Vert \nabla f( \bar \vy^{(k)}) \right\Vert \le \frac{F_0}{DKT}+ \frac{R_T\left(\vu^{(1)},\cdots,\vu^{(K)}\right)}{DKT} +  L T^2 D^2,
\end{align*}
where $F_0 = f(\vx_0) - f^*$ and $\bar \vy^{(k)} = {\sum_{n = (k-1)T+1}^{kT} \vy_n}/{T}$.
\end{lem}

\subsection{Inner Loop: Lazy Extra-Newton}

\begin{algorithm}[t]  
\caption{NALEN ($f,\vx_0, N$)}  \label{alg:ALEN-NC}
\begin{algorithmic}[1] 
\STATE Set the parameters $T \in \sN_+$, $\eta,D \in \sR_+$ according to Eq. (\ref{eq:para-ALEN-NC}). \\[1mm]
\STATE $\vv_0 = \Delta_0 =  - \frac{D \nabla f(\vx_0)}{\Vert \nabla f(\vx_0) \Vert}$, ~ $K = \lfloor \frac{N}{T} \rfloor$\\[1mm]
\STATE \textbf{for} $n=0,\cdots,N-1$\\[1mm]
\STATE ~~ $ \vz_{n} =\vx_n + \frac{1}{2} \Delta_n$ \\[1mm]
\STATE  ~~ $\Delta_{n+1} = \argmin_{\Delta \in \sB^d_{D}(\vzero)}
\langle \nabla f(\vz_n),\Delta \rangle +
\frac{1}{4}\langle \nabla^2 f(\vz_{\pi(n)}) (\Delta - \Delta_n), \Delta - \Delta_n \rangle  + \frac{1}{2 \eta} \Vert \Delta - \vv_n \Vert^2$ \label{line:TR-sub} 
~~$\triangleright$ Implementation details in Section \ref{sec:imple}.  \\[1mm] 
\STATE~~ $ \vx_{n+1} = \vx_n + \Delta_{n+1}$~ $ \vy_{n+1} = \vx_n + \frac{1}{2} \Delta_{n+1}$, ~ $\vg_{n+1} = \nabla f(\vy_{n+1})$   \\[1mm]
\STATE \label{line:extragrad} ~~$\vv_{n+1} =\argmin_{\vv \in \sB^d_D(\vzero)} \langle \vg_{n+1}, \vv  \rangle + \frac{1}{2 \eta} \Vert \vv - \vv_n \Vert^2$ \\
~~~~~~~~~~~$= {\rm Clip}_D( \vv_n - \eta \vg_{n+1})$,  ~~  where ${\rm Clip}_D(\vv) = \min \{ D / \Vert \vv \Vert, 1 \} \cdot \vv$  \\[1mm]
\STATE \textbf{end for} \\[1mm]
\STATE Set $\bar \vy^{(k)} = \frac{1}{T}  \sum_{t=1}^T \vy_{(k-1)T+t},~ \forall k = 1,\cdots,K$. \\[1mm]
\STATE \textbf{return} $\hat \vy = \arg \min_{\vy \in \{ \bar \vy^{(1)}, \cdots, \bar \vy^{(K)} \} } \Vert \nabla f(\vy) \Vert $ 
\end{algorithmic}
\end{algorithm}


We propose our algorithm for nonconvex optimization in Algorithm \ref{alg:ALEN-NC}.
Differing from previous works that solve the online learning subproblem in Section~\ref{subsec:back-o2nc} using optimistic online gradient \citep{cutkosky2023optimal} or quasi-Newton method \citep{jiang2025improved}, we apply a modified LEN algorithm \citep{chen2025secondorder} to solve this subproblem, which leads to a faster convergence rate by explicitly exploiting second-order information. Compared to LazyCRN \citep{doikov2023second} for offline setups, LEN takes an extra gradient step to obtain sublinear regret for online setups.

In fact, the LEN algorithm \citep{chen2025secondorder} is a special instance of the optimistic online gradient methods with hints \citep{rakhlin2013online,hazan2010extracting,chen2021impossible}, which takes the following extragradient update \citep{korpelevich1976extragradient}, which, given domain~$\gX$, direction variable~$\Delta$, loss vector $\vg_n$, takes the update
\begin{align} \label{eq:eg-concept}
\begin{split}
      \Delta_{n+1} =& \argmin_{\Delta \in \gX}
    \langle \vh_{n+1}, \Delta \rangle + \frac{1}{2 \eta_n} \Vert \Delta - \vv_n \Vert^2,
 \\
    \vv_{n+1} =& \argmin_{\vv \in \gX} \langle \vg_{n+1}, \vv  \rangle + \frac{1}{2 \eta_n} \Vert \vv - \vv_n \Vert^2,
\end{split}
\end{align}
where the hint vector $\vh_{n+1} \in \sR^d$ is a prediction of the future gradient $\vg_{n+1}$, and we have
 $\gX = \sB_D^d(\vzero)$ and $\vg_{n+1} = \nabla f(\vy_{n+1})$ in our setup. The following standard lemma of extragradient update shows that, if $\vh_{n+1} \approx \vg_{n+1}$, then the above updates yield a small regret.
\begin{lem} \label{lem:eg-update}
Consider the update formula in Eq. (\ref{eq:eg-concept}). For any $\vu \in \gX$, we have
\begin{align} 
\begin{split}
    \langle \vg_{n+1}, \Delta_{n+1} - \vu \rangle \le& \frac{1}{2 \eta_n} \Vert \vu - \vv_n \Vert^2 - \frac{1}{2 \eta_n} \Vert \vu - \vv_{n+1} \Vert^2 +\eta_n
     \Vert \vh_{n+1} - \vg_{n+1} \Vert^2 \\
     -& \frac{1}{2 \eta_n} \Vert \Delta_{n+1} - \vv_n \Vert^2 - \frac{1}{4 \eta_n} \Vert \Delta_{n+1} - \vv_{n+1}  \Vert^2 .
\end{split}
\end{align}
\end{lem}

Let $\vz_n = \vx_n + \Delta_n /2$ be the extrapolated point.
We design the hint vector $\vh_{n+1}$ as
\begin{align} \label{eq:hint-vec}
    \vh_{n+1} := \nabla f(\vz_n) + \nabla^2 f(\vz_{\pi(n)}) (\vy_{n+1} - \vz_n),
\end{align}
which is modified from the construction $\vh_{n+1} = \nabla f (\Delta_n) + \nabla^2 f(\Delta_{\pi(n)})(\Delta_{n+1} - \Delta_n)$ in \citep{chen2025secondorder} to reduce the prediction error $\Vert \vh_{n+1} - \nabla f(\vy_{n+1}) \Vert$ in our setup.
We can show that $\Vert \vh_{n+1} - \nabla f(\vy_{n+1}) \Vert = \gO(m L D \Vert \Delta_{n+1} - \Delta_n \Vert )$. 
Then we can substitute Lemma~\ref{lem:eg-update} into Lemma~\ref{lem:reg} and obtain the following regret bound for $R_T\left(\vu^{(1)},\cdots,\vu^{(K)} \right)$ defined in Eq. (\ref{dfn:shifting-reg}).

\begin{lem} \label{lem:reg}
Under Assumption \ref{asm:Hessian-Lip}, running Algorithm \ref{alg:ALEN-NC} with stepsize $\eta =\Theta(1/ (mLD))$
ensures 
\begin{align} \label{eq:reg-bound}
    R_T\left(\vu^{(1)},\cdots,\vu^{(K)} \right)= \gO(m K L D^3).
\end{align}
\end{lem}
Note that the extragradient updates (\ref{eq:eg-concept}) with the hint vector $\vh_{n+1}$ in Eq. (\ref{eq:hint-vec}) are implementable. The first step, \textit{i.e.}, a gradient step using the hint vector is equivalent to Line \ref{line:TR-sub},
which can be verified by using the first-order optimality condition (Proposition \ref{prop:impli-expli}). It
corresponds to solving a trust region subproblem \citep{conn2000trust}, whose implementation details are presented in Section \ref{sec:imple}. The second step, \textit{i.e.}, an (projected) extragradient step with $\vg_{n+1}$ on domain $\sB_D^d(\vzero)$ is equivalent to Line \ref{line:extragrad}, a clipped gradient step.
Finally, we obtain our NALEN algorithm by solving the online learning subproblem for the procedure described in Eq. (\ref{eq:eg-concept}) and (\ref{eq:hint-vec}), and we present the resulting method in Algorithm~\ref{alg:ALEN-NC}. 
By combining the O2NC guarantee in Lemma \ref{lem:o2nc}  and appropriately tuning the parameters $T$ and~$D$, we can obtain the following theorem for this method. 

\begin{thm} \label{thm:ALEN-NC}
Under Assumption \ref{asm:Hessian-Lip} and \ref{asm:nonconvex}, running Algorithm \ref{alg:ALEN-NC} with parameters
\begin{align} \label{eq:para-ALEN-NC}
    \eta =  \Theta\left( \frac{1}{mLD} \right), \quad T = \Theta(m^{1/3}), \quad {\rm and} \quad   D =\Theta\left( \left( \frac{F_0}{N L m^{2/3}} \right)^{1/3}  \right) 
\end{align}
finds an $\epsilon$-stationary point in $N = \gO( m + m^{1/3} F_0 L^{1/2} \epsilon^{-3/2} )$ iterations.


In particular, choosing $m = \Theta(\bar d)$ to minimize the total oracle costs in Eq. (\ref{eq:total-oracle-cost}), the equivalent gradient complexity of Algorithm \ref{alg:ALEN-NC} is $\gO( \bar d + \bar d^{1/3} F_0 L^{1/2} \epsilon^{-3/2} )$.
\end{thm}

Our result in Theorem \ref{thm:ALEN-NC} significantly improves the one of LazyCRN \citep{doikov2023second} by a factor of $\bar d^{1/6}$ and the one of RAH-AGD \citep{adil2025balancing} by a factor of $\ln^{18} \epsilon^{-1}$. In the next section, we will further build on this algorithm to propose a faster algorithm for convex optimization. 


\section{Results for Convex Optimization}

Building upon the NALEN algorithm for nonconvex problems, we further propose a double-loop method CALEN for convex problems. In the outer loop (Algorithm~\ref{alg:A-LEN}), our algorithm applies the Monteiro-Svaiter acceleration \citep{monteiro2013accelerated} that requires an MS oracle (formally introduced in Definition \ref{dfn:MS-oracle}), which can be obtained by making the gradient of a proximal objective small. Then, in the inner loop (Algorithm \ref{alg:sub-solver}), our algorithm calls NALEN to minimize the proximal objective and obtain an MS oracle efficiently. Now, let us introduce each loop in detail.


\subsection{Outer Loop: Optimal Accelerated Newton Proximal Extragradient} \label{subsec:MS-acc-intro}

We present the outer loop of our method in Algorithm \ref{alg:A-LEN}, which is the same as the recently proposed optimal A-NPE (OptANPE) algorithm \citep{carmon2022optimal}, an improved version of the original A-NPE method \citep{monteiro2013accelerated}. Now, let us give a brief introduction to both the A-NPE and OptANPE method. Roughly speaking, both of them performs the following iterations, which first update the primal variable~$\vy_t$ via a second-order proximal point step from the extrapolated point $\bar \vx_t = (A_t\vx_t + a_{t+1} \vv_t)/A_{t+1} $ in Eq. (\ref{eq:MS-update-extrapola}), (\ref{eq:MS-update-prox}), and then update the dual variable $\vv_{t+1}$ via an extragradient step in Eq. (\ref{eq:ANPE-extra-step}). In the seminal work, \citet{monteiro2013accelerated} showed that the following conceptual update can ensure a fast convergence rate of $\gO((\gamma/\epsilon)^{2/7})$.
\begin{tcolorbox}[colback=white, colframe=black, boxrule=0.2mm] 
\vspace{-4mm}
\begin{align}
    \bar \vx_t &= \frac{A_t}{A_{t} + a_{t+1}} \vx_t + \frac{a_{t+1}}{A_t +a_{t+1}} \vv_t
    \label{eq:MS-update-extrapola}, ~~ {\rm where}~~ a_{t+1} = \frac{1}{2}(\lambda_{t} + \sqrt{1 + \lambda_{t} A_t}) \\
    \vy_t &\approx \arg \min_{\vy \in \sR^d} \left\{ f(\vy) + \frac{\gamma}{3} \Vert \vy - \bar \vx_t \Vert^3  \right\} \quad {\rm and} \quad \lambda_t =  \frac{1}{\gamma\Vert \vy_t - \bar \vx_t \Vert} \label{eq:MS-update-prox}\\
    \vv_{t+1} &= \vv_t - a_{t+1} \nabla f(\vy_t), \quad A_{t+1} = A_t + a_{t+1}, \quad {\rm and} \quad \vx_{t+1} = \vy_t \label{eq:ANPE-extra-step}
\end{align}
\end{tcolorbox}

There are two challenges in implementing the above updates. One subtlety is that the proximal point step in Eq. (\ref{eq:MS-update-prox}) usually can not be exactly solved, and we can only obtain an inexact solution of it. To address this issue, \citet{monteiro2013accelerated} introduced the following error condition such that the A-NPE method can converge. For convenience, we call this condition the MS condition, and the oracle that can output an acceptable solution that satisfies the MS condition the MS oracle.

\begin{dfn}[MS oracle for \eqref{eq:MS-update-prox}] \label{dfn:MS-oracle}
Given parameters $\sigma \in (0,1)$ and $\gamma>0$, we define
$\sO^{\rm MS}_{f,\sigma,\gamma}: \sR^d \rightarrow \sR^d$ as an MS oracle for function $f: \sR^d \rightarrow \sR$, which takes a query point
$\vx \in \sR^d$ and returns an approximate proximal point $\hat \vy = \sO^{\rm MS}_{f,\sigma,\gamma}(\vx)$ that satisfies the following MS condition: 
\begin{align} \label{eq:MS-oracle}
\Vert \nabla f_{\vx, \gamma}(\hat \vy) \Vert \le  \sigma \gamma \Vert \vx -\hat \vy \Vert^2,
\end{align}
where $f_{\vx, \gamma}: \sR^d \rightarrow \sR^d$ is the second-order proximal function defined as
\begin{align} \label{eq:second-order-prox}
    f_{\vx, \gamma}(\vy) :=  f(\vy) + \gamma \Vert \vy - \vx \Vert^{3} /3.
\end{align}
\end{dfn}

\begin{algorithm*}[t]  
\caption{CALEN$(f, \vx_0, T) $
}\label{alg:A-LEN}
\begin{algorithmic}[1] 
\STATE Set the parameters $\sigma \in (0,1)$, $\gamma \in \sR_+$, $S,N \in \sN_+$ according to Theorem \ref{thm:ALEN}.
\STATE $\vv_0 = \vx_0$, $\bar \vx_0 = \vx_0$, $A_0 = 0$ \\
\STATE Run Algorithm \ref{alg:sub-solver} with parameters $S,N$ to obtain $ \vy_0 = {\sO^{\rm MS}_{f,\sigma,\gamma}(\vx_0)}$. \label{line:MS-ini}  \\
\STATE Set $\lambda_0' = \lambda_0= 1/(\gamma \Vert \vy_0 - \vx_0 \Vert) $. \\
\FOR{$t = 0,\cdots, T-1$}  
\STATE \quad $a_{t+1}' =  ( \lambda_{t}' + \sqrt{1 + 4 \lambda_{t}' A_t}) / 2$, ~ $A_{t+1}' = A_t + a_{t+1}'$,~ $\bar \vx_t = (A_t \vx_t + a_{t+1}' \vv_t)/{A_{t+1}'} $ \label{line:momentum} \\
\STATE \quad \textbf{if} $t > 0 $ \textbf{then} 
\STATE \quad \quad Run Algorithm \ref{alg:sub-solver} with parameters $S,N$ to obtain $\vy_t = {\sO^{\rm MS}_{f, \sigma, \gamma}(\bar \vx_t)}$. \label{line:MS} \\
\STATE \quad \quad Set $\lambda_t = 1/(\gamma \Vert \vy_t - \bar \vx_t \Vert)$. \\
\STATE \quad \textbf{end if} \\
\STATE \quad  \textbf{if} $\lambda_{t} \ge  \lambda_{t}'$ \textbf{then} \label{line:begin-seacg-lambda} \\
\STATE \quad \quad $a_{t+1} = a_{t+1}' $, $A_{t+1} = A_{t+1}'$, ~ $\vx_{t} = \vy_{t} $, ~ $\lambda_{t+1}' = 2 \lambda_{t}' $\\
\STATE \quad \textbf{else} \\
\STATE \quad \quad $\beta_{t} = \lambda_{t}/\lambda_{t}'$,~ $a_{t+1} = \beta_{t} a_{t+1}'$, ~ $A_{t+1} = A_t + a_{t+1}$ \\
\STATE \quad \quad $ \vx_{t+1} = ((1- \beta_{t}) A_t \vx_t + \beta_{t} A_{t+1}' \vy_{t}) / A_{t+1}$,~ $\lambda_{t+1}' = \lambda_{t}'/2$ \\
\STATE \quad \textbf{end if} \label{line:end-search-lambda} \\
\STATE \quad $\vv_{t+1} = \vv_t - a_{t+1} \nabla f(\vy_{t})$ \label{line:extra-ALEN} \\
\ENDFOR \\
\RETURN $\vx_T$ 
\end{algorithmic}
\end{algorithm*}

Another subtlety of the A-NPE method is that 
Eq. (\ref{eq:MS-update-extrapola}) and (\ref{eq:MS-update-prox}) jointly  define an implicit equation for $(\vy_t, \lambda_{t})$, since $ \vy_t$ is calculated with $\bar \vx_t$ in Eq. (\ref{eq:MS-update-prox}), which depends on $\lambda_t$ in Eq. (\ref{eq:MS-update-extrapola}), which in turn depends on $ \vy_t$ and $\bar \vx_t$ through $\lambda_t =  1/(\gamma\Vert \vy_t - \bar \vx_t \Vert)$ in Eq.~(\ref{eq:MS-update-prox}). To address this issue and obtain an explicit method, \citet{monteiro2013accelerated} proposed to solve the implicit equation in Eq. (\ref{eq:MS-update-extrapola}) and (\ref{eq:MS-update-prox}) 
via a binary search on the variable $\lambda_t$ to find the root of the univariate function $\varphi(\lambda_t):= \lambda_t - 1/ (\gamma \Vert \vy_t - \bar \vx_t \Vert)$. Note that $\varphi(\lambda_t)$ is a continuous function in $\lambda_t$: In A-NPE updates, if $\lambda_t = 0$, then $a_{t+1} =0$ and $\bar \vx_t = \vx_t$, which implies $\varphi(\lambda_t) < 0$; if  
$\lambda_t \rightarrow + \infty$, then $a_{t+1} \rightarrow +\infty$ and $\bar \vx_t = \vv_t$, which implies $\varphi(\lambda_t) \ge 0$. Therefore, the root of this unitary function always exists in~$\sR_+$, and \citet{monteiro2013accelerated} showed that binary search terminates in $\gO(\ln \epsilon^{-1})$ iterations. \citet{carmon2022optimal} took a further step forward by conducting the binary search along with the updates. In their improved implementation of A-NPE, the algorithm adaptively adjusts a guess to $\lambda_t$ and damps the momentum when the guess overestimates. In our method, we also adopt this improved technique in Line~\ref{line:begin-seacg-lambda} to \ref{line:end-search-lambda} since it can shave off the $\ln \epsilon^{-1}$ factor in the binary search. Note that Algorithm \ref{alg:A-LEN} is exactly the same as the OptANPE method \citep{carmon2022optimal} except for the MS oracle implementation. Then we can know from \citep{carmon2022optimal} that, given access to an MS oracle, this algorithm has the following guarantee.



\begin{lem}[{\citet[Theorem 1 with $s=\alpha=2$]{carmon2022optimal}}] \label{lem:A-NPE}
Under Assumption \ref{asm:convex}, if Line~\ref{line:MS-ini} and~\ref{line:MS} always return an MS oracle $\sO^{\rm MS}_{f,\sigma,\gamma}$ with $\sigma \in (0, 0.99]$ then running Algorithm \ref{alg:A-LEN} guarantees $f(\vx_T) - f(\vx^*) \le \epsilon$ in $T = \mathcal{O}(( \gamma R_0^{3} / \epsilon )^{{2}/{7}})$ iterations, where 
$R_0 = \Vert \vx_0 - \vx^* \Vert$ is the initial distance to an arbitrary minimizer $\vx^*$.
\end{lem}


\subsection{Inner Loop: Nonconvex Accelerated Lazy Extra-Newton} \label{subsec:ALEN-algo}


\begin{algorithm*}[t]  
\caption{$\sO^{\rm MS}_{f,\sigma,\gamma}(\vx)$ via NALEN on $f_{\vx,\gamma}$}  \label{alg:sub-solver}
\begin{algorithmic}[1] 
\renewcommand{\algorithmicrequire}{ \textbf{Parameters:}}
\STATE Set the parameters $S$ and $N$ according to Eq. (\ref{eq:para-sub-solver}).  \\ [1mm]
\STATE $\vy_0 = \arg \min_{\vy \in \sR^d} \langle \nabla f(\vx), \vy \rangle + \langle \nabla^2 f(\vx) (\vy - \vx), \vy-\vx \rangle/2 + (L + 2\gamma) \Vert \vy - \vx \Vert^3 /3$ \label{lin:ini-CRN} \\[1mm] 
$\triangleright$ Implementation details in Section \ref{sec:imple}.  \\[1mm] 
\STATE \textbf{for} $s = 0,1,\cdots, S$ \\ [1mm]
\STATE \quad $\vy_{s+1} = \text{NALEN}(f_{\vx,\gamma}, \vy_s, N)$ \\[1mm]
\STATE \textbf{end for} \\[1mm]
\STATE \textbf{return} $\hat \vy = \vy_{S+1}$  
\end{algorithmic}
\end{algorithm*}

We present our implementation of MS oracle in Algorithm~\ref{alg:sub-solver}.
Differing from prior works \citep{carmon2022optimal,nesterov2021inexact,nesterov2023inexact}
that implement the MS oracle $\sO_{f,\sigma,\gamma}^{\rm MS}(\vz)$ using a single CRN step, we apply a restart variant of the NALEN method (Algorithm \ref{alg:ALEN-NC}) after the CRN step in Line \ref{lin:ini-CRN} of Algorithm \ref{alg:sub-solver}. Using lazy Hessian updates, we can ensure that the additional NALEN iterations remain comparable costs to the CRN step by choosing the parameters $S = \Theta(\ln m)$ and $N = \Theta(m)$, but returns a more aggressive MS oracle $\sO_{f,\sigma,\gamma}^{\rm MS}(\vz)$ with $\gamma = \Theta(  L/m^{4/3} )$, improving the $\gamma = \Theta(L)$ in previous implementations \citep[Section~3.1]{carmon2022optimal}. We find our idea very similar to the Catalyst acceleration \citep{lin2018catalyst} in first-order methods, which has been used in many works \citep{frostig2015regularizing,ivanova2021adaptive,carmon2022recapp}.
However, to the best of our knowledge, this idea has not been fully explored in second-order optimization. 
 
Now, let us start our analysis of CALEN by presenting some useful properties of $f_{\vx,\gamma}(\vy)$.

\begin{lem} \label{lem:sub-prob-grad-dom}
Under Assumption \ref{asm:Hessian-Lip} and \ref{asm:convex}, the function $f_{\vx,\gamma}(\vy)$ defined in Eq. (\ref{eq:second-order-prox}) satisfies:
\begin{enumerate}
    \item  $f_{\vx,\gamma}(\vy)$ has $(L+2 \gamma)$-Lipschitz continuous Hessians.
    \item Let $\vy_\gamma^*(\vx) = \arg \min_{\vy \in \sR^d} f_{\vx,\gamma}(\vy)$ and $f_{\vx,\gamma}^* = \min_{\vy \in \sR^d} f_{\vx,\gamma}(\vy)$. For any $\vy \in \sR^d$,  
    \begin{align}
        \Vert \nabla f_{\vx, \gamma}(\vy) \Vert^{3/2} \ge \max\left\{ \frac{3}{2} \sqrt{\frac{\gamma}{2}} (f_{\vx, \gamma}(\vy) - f_{\vx,\gamma}^*), \left( 
        \frac{\gamma}{2} \Vert \vy - \vy_\gamma^*(\vx) \Vert^2
        \right)^{3/2} \right\}. \label{eq:cubic-grad-uniform}
    \end{align}
\end{enumerate}
\end{lem}

By exploiting the above properties of $f_{\vx,\gamma}(\vy)$  stated in Lemma \ref{lem:sub-prob-grad-dom} and invoking Theorem \ref{thm:ALEN-NC} for NALEN, we can show that our sub-solver (Algorithm \ref{alg:sub-solver}) can implement an MS oracle efficiently.


\begin{lem} \label{lem:sub-solover}
Under Assumption \ref{asm:Hessian-Lip} and \ref{asm:convex}, given any $\sigma,\gamma>0$ and $m \ge 2$, running Algorithm \ref{alg:sub-solver} with parameters
\begin{align} \label{eq:para-sub-solver}
    S = \Theta\left(\ln \left( \frac{1}{\sigma} \left( \frac{L+2\gamma}{\gamma} \right)^{2/3} \right) \right ) \quad \text{and} \quad N= \Theta \left( m^{1/3} \sqrt{\frac{L+2\gamma}{\gamma}} \right)
\end{align}
can implement an MS oracle $\sO_{f,\sigma,\gamma}^{\rm MS}(\vz)$ in $SN$ total iterations.
\end{lem}

According to Lemma \ref{lem:sub-solover}, for any $m \ge 2$, the optimal Hessian-gradient trade-off in the MS-Solver is achieved under the setting $\gamma = \Theta( L / m^{4/3})$. In this case, Algorithm \ref{alg:sub-solver} can return an MS oracle $\sO_{f,\sigma,\gamma}^{\rm MS}$ using $SN = \gO(m \ln m)$ iterations. Moreover, invoking Theorem \ref{thm:ALEN} with $\gamma = \Theta( L / m^{4/3})$, we know that Algorithm \ref{alg:A-LEN} finds an $\epsilon$-solution of Problem~(\ref{prob:main}) in $T = \mathcal{O}( m^{-8/21} \epsilon^{-2/7} )$ iterations, and the total iteration complexity of Algorithm \ref{alg:A-LEN} is $TSN = \gO( (m + m^{13/21} \epsilon^{-2/7}) \ln m )$. We summarize the result in the following theorem.


\begin{thm} \label{thm:ALEN}
Under Assumption \ref{asm:Hessian-Lip} and \ref{asm:convex}, setting the parameters $\gamma = \Theta(L / m^{4/3})$, $\sigma \in (0, 0.99]$, and $S, N$ according to Eq. (\ref{eq:para-sub-solver}), then Algorithm \ref{alg:A-LEN} can find an $\epsilon$-global solution in total $TSN = \gO( (m + m^{13/21} (L R_0^3/\epsilon)^{2/7}) \ln m )$ iterations.

In particular, choosing $m = \Theta(\bar d)$ to minimize the total oracle costs in Eq. (\ref{eq:total-oracle-cost}), the equivalent gradient complexity of Algorithm \ref{alg:A-LEN} is $ \gO( (\bar d + \bar d^{13/21} (L R_0^3/\epsilon)^{2/7}) \ln \bar d)$.
\end{thm}

Our result in Theorem \ref{thm:ALEN}  improves the $\gO\left(\bar d^{1/2} (L R_0^3 / \epsilon)^{1/2}\right)$ of LazyCRN \citep{chayti2023unified} when $d = \gO( \epsilon^{-1})$, and matches the $\gO(\bar d)$ of LazyCRN up to only an extra $\ln \bar d$ factor when $d = \Omega(\epsilon^{-1})$. It is open whether one can shave off this additional $\ln \bar d$ factor, while still achieving the optimal $\epsilon^{-2/7}$ dependency in $\epsilon$ \citep{arjevani2019oracle} with lazy Hessian updates. Additionally, we extend our algorithm to strongly convex problems using restart mechanisms in Appendix~\ref{apx:ALEN-SC}.

\section{Implementation Details} \label{sec:imple}

This section provides implementation details of Newton subproblems in our algorithms. In Line \ref{line:TR-sub} of Algorithm \ref{alg:ALEN-NC}, we need to solve the 
following trust-region Newton subproblem \citep{conn2000trust}:
\begin{align} \label{eq:tr-sub-prob}
    \vh^\ast = \arg \min_{\vh \in \sB_D^d(\vzero)}\langle \vb, \vh \rangle  + \langle \mA \vh, \vh \rangle /2.
\end{align}
In Line \ref{lin:ini-CRN} of Algorithm \ref{alg:sub-solver}, we need to solve the CRN subproblem \citep{nesterov2006cubic}:
\begin{align}  \label{eq:CRN-prob}
    \vh^* = \arg \min_{\vh \in \sR^d} \langle \vb, \vh \rangle +\langle \mA \vh, \vh \rangle /2 + M \Vert \vh \Vert^3 /6.
\end{align}
It is known that both the global solutions $\vh^*$ of these subproblems satisfy the following fact.
\begin{lem} \label{lem:sub-problem}
Let $\mA \in \sR^{d \times d}$ be symmetric.
For either Problem (\ref{eq:tr-sub-prob}) or (\ref{eq:CRN-prob}), there always
exists a continuous, monotone, convex function $\varphi(\tau): (\underline{\tau}, +\infty) \rightarrow \sR$ with $\underline{\tau} = \max\{ - \lambda_{\min} (\mA),0 \}$, such that its root  $\tau^* \in \sR_+$ satisfies 
\begin{align} \label{eq:root-cond}
   (\mA + \tau^* \mI_d) \vh^* = -\vb  \quad {\rm and} \quad  \mA +\tau^* \mI_d \succeq \mO_d.
\end{align}
\end{lem} 

The root $\tau^*$ in the above can be efficiently solved by the univariate Newton method \citep[Appendix A.1]{nesterov2018lectures} or a standard binary search procedure. Moreover, for convex problems where $\mA \succeq \mO_d$, one can also apply a binary search on the logarithmic space \citep{monteiro2012iteration} for faster convergence. After the calculation of $\tau^*$ and omitting the degenerate cases that can be avoided by imposing small perturbations \citep{nesterov2006cubic}, we can then solve these subproblems as a normal Newton step $\vh^* =  -(\mA + \tau^* \mI_d)^{-1} \vb$. 
Moreover, as suggested by \citet[Section 6.1]{doikov2023second}, when using lazy Hessian updates, we only need to calculate a spectral decomposition of $\mA$ at the snapshot point, then the matrix inversion in each step can be implemented in $\gO(d^2)$ arithmetic operations. Note that the spectral decomposition can be implemented in matrix multiplication time $\gO(d^{\omega})$ \citep{demmel2007fast}, where the current best exponent $\omega$ is around $2.37$ \citep{duan2023faster}. Hence, if we further add the extra costs of solving these Newton subproblems into Eq.~(\ref{eq:total-oracle-cost}), the total computational complexity is $\gO\left( n_\mH(c_\mH + d^\omega ) + n_\vg (c_\vg + d^2) \right)$, indicating that the optimal choice of $m$ is $m = \Theta( (c_\mH + d^\omega ) / (c_\vg+ d^2)) $. 




\section{Conclusion and Future Work}

In this work, we propose NALEN and CALEN algorithms and show that they achieve the equivalent gradient complexity of $\gO( \bar d + \bar d^{1/3} \epsilon^{-3/2})$ and  $\gO( (\bar d + \bar d^{13/21} \epsilon^{-2/7}) \ln \bar d)$ for nonconvex and convex optimization, where $\bar d \le d$ is the effective dimension of the problem.

It is open whether the dependencies on $\bar d$ we achieve are optimal, and it is important to seek improved upper
and lower bounds in the future. It is also interesting to reduce the computational costs of second-order methods for other problem classes, such as 
self-concordant problems \citep{nesterov1994interior},
third-order smooth problems \citep{nesterov2021implementable,nesterov2021superfast,nesterov2021inexact,gasnikov2019near}, or Hessian-stable problems \citep{karimireddy2018global,carmon2020acceleration}.

\section*{Acknowledgment}
Jingzhao Zhang is supported by National Key R\&D Program of China 2024YFA1015800 and the Shanghai Qi Zhi Institute Innovation Program. 
Luo Luo is supported by the National Natural Science Foundation of China (No.12571557),
the Major Key Project of PCL under Grant PCL2024A06, and the Shanghai Basic Research Program (23JC1401000). Chengchang Liu is supported by the National Natural Science Foundation of China (624B2125). 

\bibliography{sample}

\newpage 

\appendix

\section{Proof of Lemma \ref{lem:eg-update}}

\begin{proof}
By the first-order optimality condition in the 
two steps 
in the update formula (\ref{eq:eg-concept}),we have
\begin{align*}
      0 \le& \langle \eta_n \vh_{n+1} + \Delta_{n+1}- \vv_n, \vu' - \Delta_{n+1} \rangle, \quad \forall \vu' \in \gX; \\
    0 \le& \langle \eta_n \vg_{n+1} + \vv_{n+1}- \vv_n, \vu - \vv_{n+1} \rangle, \quad \forall \vu \in \gX.
\end{align*}
Let $\vu' = \vv_{n+1}$. Then we can use the two inequalities above to obtain
\begin{align*}
    & \langle \vg_{n+1}, \Delta_{n+1} - \vu \rangle \\
    = &   \langle \vg_{n+1}, \vv_{n+1} - \vu \rangle  + \langle \vg_{n+1}, \Delta_{n+1} - \vv_{n+1} \rangle \\
    = & \langle \vg_{n+1}, \vv_{n+1} - \vu \rangle  +  \langle \vh_{n+1}, \Delta_{n+1} - \vv_{n+1} \rangle \\
    & + \langle  \vg_{n+1} - \vh_{n+1}, \Delta_{n+1} - \vv_{n+1}  \rangle \\
    \le & \frac{1}{\eta_n} \langle \vv_{n+1} - \vv_n , \vu - \vv_{n+1} \rangle + 
    \frac{1}{\eta_n}
    \langle\Delta_{n+1} - \vv_n , \vv_{n+1} - \Delta_{n+1} \rangle \\
    &+  \langle  \vg_{n+1} - \vh_{n+1}, \Delta_{n+1} - \vv_{n+1}  \rangle.
\end{align*}
Next, using the identity $ \langle \va, \vb \rangle =  \frac{1}{2} \Vert \va+\vb \Vert^2- \frac{1}{2} \Vert \va \Vert^2 - \frac{1}{2} \Vert \vb \Vert^2$ on the first two terms and the Young's inequality $\langle \va,\vb \rangle \le \gamma \Vert \va \Vert^2 + \frac{1}{4 \gamma} \Vert \vb \Vert^2 $ for any $\gamma>0$ on the third term, we further have 
\begin{align*}
     & \langle \vg_{n+1}, \Delta_{n+1} - \vu \rangle \\ 
     \le & \frac{1}{2 \eta_n} \Vert \vu - \vv_n \Vert^2 -  \bcancel{ \frac{1}{2 \eta_n} \Vert \vv_{n+1} - \vv_n \Vert^2} - \frac{1}{2 \eta_n} \Vert \vu - \vv_{n+1} \Vert^2 \\
     & + \bcancel{\frac{1}{2 \eta_n} \Vert \vv_{n+1} - \vv_n \Vert^2} - \frac{1}{2 \eta_n} \Vert \Delta_{n+1} - \vv_n \Vert^2 - \frac{1}{2 \eta_n} \Vert \Delta_{n+1}  -\vv_{n+1}  \Vert^2 \\
     & + \eta_n
     \Vert \vg_{n+1} - \vh_{n+1} \Vert^2 + \frac{1}{4 \eta_n} \Vert \Delta_{n+1} - \vv_{n+1}  \Vert^2 \\
     =& \frac{1}{2 \eta_n} \Vert \vu - \vv_{n} \Vert^2 - \frac{1}{2 \eta_n} \Vert \vu - \vv_{n+1} \Vert^2 \\
     &  - \frac{1}{2 \eta_n} \Vert \Delta_{n+1} - \vv_n \Vert^2 - \frac{1}{4 \eta_n} \Vert \Delta_{n+1} -\vv_{n+1}  \Vert^2  + \eta_n
     \Vert \vh_{n+1} - \vg_{n+1} \Vert^2.
\end{align*}
\end{proof}

\section{Proof of Lemma \ref{lem:reg}} \label{apx:proof-reg}

\begin{proof} Under Assumption \ref{asm:Hessian-Lip}, we have
\begin{align} \label{eq:Hess-lip}
        \Vert \nabla^2 f(\vx) - \nabla^2 f(\vx') \Vert \le L \Vert \vx-  \vx' \Vert, \quad \forall \vx,\vx' \in \sR^d.
    \end{align}
Then, a consequence of Taylor expansion is the following fact \citep[Lemma 1.2.4]{nesterov2018lectures}.
\begin{align}  \label{eq:Hessian-Lip-taylor}
    \Vert \nabla f(\vx') - \nabla f(\vx) - \nabla^2 f(\vx) (\vx' - \vx) \Vert \le \frac{L}{2} \Vert \vx - \vx' \Vert^2, \quad \forall \vx, \vx' \in \sR^d.
\end{align}
Now, note that Algorithm \ref{alg:ALEN-NC} is equivalent to using the optimistic online gradient descent described in Eq.~(\ref{eq:eg-concept}) to solve the online learning sub-problem described in Section \ref{subsec:back-o2nc} with 
constant stepsize~$\eta_n = \eta$ and hint vector $\vh_{n+1}  =  \nabla f(\vz_n) + \nabla^2 f(\vz_{\pi(n)}) (\vy_{n+1} - \vz_n)$. We want to use Lemma \ref{lem:eg-update} to prove the regret bound, which requires an upper-bound of the prediction error $\Vert \vh_{n+1} - \nabla f(\vy_{n+1}) \Vert$. We analyze it by starting from the following identity
\begin{align*}
    &\vh_{n+1} -\nabla f(\vy_{n+1})\\
    =&  \nabla f(\vz_n ) + \nabla^2 f({\vz_{\pi(n)}}) ( \vy_{n+1} - \vz_n) - \nabla f(\vy_{n+1}) \\
    =& \nabla f(\vz_n ) + \nabla^2 f({\vz_{n}}) ( \vy_{n+1} - \vz_n) - \nabla f(\vy_{n+1})  + (\nabla^2 f({\vz_{\pi(n)}}) -\nabla^2 f({\vz_{n}})) ( \vy_{n+1} - \vz_n),
\end{align*}
where the second equality uses the relationship $\vy_{n+1} - \vz_n =  \frac{1}{2} (\Delta_{n+1} - \Delta_n)$. Then,
taking the norm on both sides leads to
\begin{align*}
    & \Vert \vh_{n+1} -\nabla f(\vy_{n+1}) \Vert \\
    \le& \Vert \nabla f(\vz_n ) + \nabla^2 f({\vz_{n}}) ( \vy_{n+1} - \vz_n) - \nabla f(\vy_{n+1})   \Vert + \Vert (\nabla^2 f({\vz_{\pi(n)}}) -\nabla^2 f({\vz_{n}})) ( \vy_{n+1} - \vz_n) \Vert \\
    \le& \frac{L}{2} \Vert  \vy_{n+1} - \vz_n \Vert^2 + L \Vert \vz_{\pi(n)} - \vz_n \Vert \Vert \vy_{n+1} - \vz_n \Vert \\
    =& \frac{L}{4} \Vert \vy_{n+1} - \vz_n \Vert \Vert \Delta_{n+1} - \Delta_n \Vert + \frac{L}{2} \Vert \vz_{\pi(n)} - \vz_n \Vert  \Vert \Delta_{n+1} - \Delta_n \Vert,
\end{align*}
where the first inequality uses the triangle inequality $\Vert \va + \vb \Vert \le \Vert \va \Vert + \Vert \vb \Vert$, the second one uses inequalities (\ref{eq:Hess-lip}) and (\ref{eq:Hessian-Lip-taylor}), the third inequality uses $\vy_{n+1} - \vz_n =  \frac{1}{2} (\Delta_{n+1} - \Delta_n)$ again. Next, recalling $\vz_n -\vz_{\pi(n)} = \sum_{i = \pi(n)}^{n-1} \Delta_i + \frac{1}{2} (\Delta_n - \Delta_{\pi(n)}) $ and $\Vert \Delta_n \Vert \le D$ for all $n$, we have
\begin{align*}
    \Vert \vy_{n+1} - \vz_n \Vert \le D \quad {\rm and} \quad \Vert \vz_n -\vz_{\pi(n)} \Vert \le m D.
\end{align*}
Therefore, we get
\begin{align*}
    \Vert \vh_{n+1} -\nabla f(\vy_{n+1}) \Vert \le& \frac{(m+1) L D}{2} \Vert \Delta_{n+1} - \Delta_n \Vert.
\end{align*}
Next, taking the square on both sides and using Young's inequality gives 
\begin{align*}
    \Vert \vh_{n+1} -\nabla f(\vy_{n+1}) \Vert^2 \le& \frac{(m+1)^2 L^2 D^2}{4} \Vert \Delta_{n+1} - \Delta_n \Vert^2 \\
    \le& \frac{(m+1)^2 L^2 D^2}{2} ( \Vert \Delta_{n+1} - \vv_n \Vert^2 + \Vert \vv_n - \Delta_n \Vert^2).
\end{align*}
Substituting the above error bound into Lemma \ref{lem:eg-update} and telescoping for $n  = (k-1)T, \cdots, kT-1$, we have
\begin{align*}
    &\sum_{n= (k-1)T}^{kT-1} \langle \nabla f(\vy_{n+1}), \Delta_{n+1} - \vu^{(k)} \rangle\\
    &\le \frac{1}{2\eta} \Vert \vv_{(k-1)T} - \vu^{(k)} \Vert^2 - \sum_{n=(k-1)T}^{kT-1} \left( 
    \frac{1}{2 \eta} - \frac{\eta (m+1)^2 L^2 D^2}{2} 
    \right) \Vert \Delta_{n+1} - \vv_n \Vert^2\\
    &\quad - \sum_{n=(k-1)T}^{kT-1} \left( 
    \frac{1}{4 \eta} -  \frac{\eta (m+1)^2 L^2 D^2}{2} 
    \right) \Vert \Delta_{n+1} - \vv_{n+1} \Vert^2 \\
    &\quad + \frac{\eta (m+1)^2 L^2 D^2}{2}  \Vert \Delta_{(k-1)T} - \vv_{(k-1)T}\Vert^2.
\end{align*}
Note that we have  $\frac{1}{4 \eta} - \frac{\eta (m+1)^2 L^2 D^2}{2}  \ge 0$ by setting $\eta \le \frac{1}{2 (m+1) LD}$  and that $\Vert \Delta_{(k-1)T}- \vv_{(k-1)T} \Vert^2 \le 2\Vert \Delta_{(k-1)T}\Vert + 2\Vert \vv_{(k-1)T} \Vert^2 \le 4D^2$; $ \Vert \vv_{(k-1)T} - \vu^{(k)}\Vert^2 \le 2\Vert \vv_{(k-1)T} \Vert^2 + 2\Vert \vu^{(k)} \Vert^2 \le 4D^2$ due to the Young's inequality.
Therefore, we can simplify the above inequality to
\begin{align*}
    \sum_{n= (k-1)T}^{kT-1} \langle \nabla f(\vy_{n+1}), \Delta_{n+1} - \vu^{(k)} \rangle \le \frac{2D^2}{\eta} + (m+1) L D^3 = \gO( m LD^3 ),
\end{align*}
where the last inequality holds because $\eta = \Theta( 1/(mLD) )$. Now, we have obtained the upper bound on the static regret in each epoch $k$. We can then easily sum up all the static regret of all epochs $k = 1,\cdots,K$ to prove that the $K$-shifting regret in Eq. (\ref{dfn:shifting-reg}) can be bounded by
\begin{align*}
    R_T\left(\vu^{(1)},\cdots,\vu^{(K)} \right)= \gO(m K L D^3)
\end{align*}
\end{proof}

\section{From Implicit to Explicit} \label{apx:impli-expli}

Eq. (\ref{eq:eg-concept}) and (\ref{eq:hint-vec}) describe an implicit algorithm, and we have shown in Lemma \ref{lem:reg} that it ensures a low regret with $R_T\left(\vu^{(1)},\cdots,\vu^{(K)} \right)= \gO(m K L D^3)$. We show in the following proposition that it actually implied an explicit algorithm that can be easily implemented in practice, which is exactly Algorithm \ref{alg:ALEN-NC} we presented in the main text.

\begin{prop} \label{prop:impli-expli}
For any convex compact set $\gX \subseteq \sR^d$,
the implicit update
\begin{align} \label{eq:impli-recall}
\begin{split}
    &\Delta_{n+1} = \argmin_{\Delta \in \gX}
    \langle \vh_{n+1}, \Delta \rangle + \frac{1}{2 \eta_n} \Vert \Delta - \vv_n \Vert^2; \\
    {\rm where} \quad &\vh_{n+1} := \nabla f(\vz_n) + \nabla^2 f(\vz_{\pi(n)}) (\vy_{n+1} - \vz_n),  
\end{split}
\end{align}
is equivalent to the following explicit update:
\begin{align} \label{eq:expli-recall}
    \Delta_{n+1} = \argmin_{\Delta \in \gX}
    \langle \nabla f(\vz_n),\Delta \rangle + \frac{1}{4}\langle \nabla^2 f(\vz_{\pi(n)}) (\Delta - \Delta_n), \Delta - \Delta_n \rangle  + \frac{1}{2 \eta_n} \Vert \Delta - \vv_n \Vert^2.
\end{align}
\end{prop}

\begin{proof}
The first-order optimality condition of the first line of Eq. (\ref{eq:impli-recall}) is
\begin{align*}
     0 \le \langle \eta_n \vh_{n+1} + \Delta_{n+1}- \vv_n, \vu - \Delta_{n+1} \rangle, \quad \forall \vu \in \gX.
\end{align*}
Next, we plug in the construction of $\vh_{n+1}$ in the second-line of Eq. (\ref{eq:impli-recall}) to obtain
\begin{align*}
    0 \le \langle \eta_n \left(\nabla f(\vz_n) + \nabla^2 f(\vz_{\pi(n)}) (\vy_{n+1} - \vz_n) \right)+ \Delta_{n+1}- \vv_n, \vu - \Delta_{n+1} \rangle, \quad \forall \vu \in \gX.
\end{align*}
Note that the variables $\vy_{n+1}$, $\vz_n$, $\Delta_{n+1}$, and $\Delta_n$ satisfy the relationship $\vy_{n+1} - \vz_n = \frac{1}{2} (\Delta_{n+1} - \Delta_n)$ in our setting. The above inequality is equivalent to
\begin{align*}
    0 \le \left \langle \eta_n\left(\nabla f(\vz_n) + \nabla^2 f(\vz_{\pi(n)}) (\Delta_{n+1} - \Delta_n) /2\right)+ \Delta_{n+1}- \vv_n, \vu - \Delta_{n+1} \right \rangle, \quad \forall \vu \in \gX,
\end{align*}
which is exactly the same as the first-order optimality condition of Eq. (\ref{eq:expli-recall}).
\end{proof}

\section{Proof of Theorem \ref{thm:ALEN-NC}}

\begin{proof}
Plugging the regret bound of Lemma \ref{lem:reg} into Lemma \ref{lem:o2nc}, we have
    \begin{align*}
     \Vert \nabla f(\hat \vy) \Vert \le \frac{1}{K} \sum_{k=1}^K \left \Vert \nabla f( \bar \vy^{(k)}) \right\Vert = \gO \left( \frac{F_0}{D N}+ \frac{m L D^2}{T} +  L T^2 D^2 \right),
    \end{align*}
We pick $T = \Theta(m^{1/3})$ to balance the second and third terms above, which leads to
 \begin{align*}
     \Vert \nabla f(\hat \vy) \Vert = \gO \left( \frac{F_0}{D N}+ m^{2/3} L D^2 \right).
    \end{align*}
Next, we can further pick $D =\Theta\left( \left( \frac{F_0}{N L m^{2/3}} \right)^{1/3}  \right)$ to balance the remaining two terms, which then leads to the final bound of
\begin{align} \label{eq:final-bound-grad}
    \Vert \nabla f(\hat \vy) \Vert = \gO \left( \left( \frac{m^{1/3}  L^{1/2} F_0}{N} \right)^{2/3} \right).
\end{align}
In other words, Algorithm \ref{alg:ALEN-NC} returns a point $\hat \vy \in \sR^d$ satisfying $\Vert \nabla f(\hat \vy) \Vert \le \epsilon$ in the gradient complexity of $\gO( m + F_0 m^{1/3} L^{1/2} \epsilon^{-3/2} )$ as the theorem claimed.
\end{proof}

\section{Proof of Lemma \ref{lem:sub-prob-grad-dom}}

\begin{proof} We define $d(\vy) = \Vert \vy \Vert^3 / 3$ as an auxiliary function and then use the properties of $d(\vy)$ to induce those of $f_{\vx,\gamma}(\vy) = f(\vy) + d(\vy-\vx)$.
\begin{enumerate}
    \item It can be verified that it has 2-Lipschitz continuous Hessians \citep[Lemma 4.2.4]{nesterov2018lectures}), which proves our first claim that $f_{\vx,\gamma}(\vy)$ has $(L+2 \gamma)$-Lipschitz continuous Hessians.
    \item   In addition, it can also be verified \citep[Lemma 4.2.3]{nesterov2018lectures} that
    \begin{align*}
    \langle \nabla d(\vy) - \nabla d(\vy') , \vy - \vy' \rangle \ge \frac{1}{2} \Vert \vy - \vy' \Vert^3, \quad \forall \vy , \vy' \in \sR^d.
\end{align*}
Since $f$ is a convex function, we have \citep[(2.1.12)]{nesterov2018lectures} that
\begin{align*}
    \langle \nabla f(\vy) - \nabla f(\vy') , \vy - \vy' \rangle \ge 0, \quad \forall \vy , \vy' \in \sR^d.
\end{align*}
Then, combining the two inequalities above gives
\begin{align}  \label{eq:grad-dom-last}
    \langle \nabla f_{\vx, \gamma}(\vy) - \nabla f_{\vx, \gamma}(\vy'), \vy - \vy' \rangle \ge \frac{\gamma}{2} \Vert \vy - \vy' \Vert^3, \quad \forall \vy , \vy' \in \sR^d.
\end{align}
Next, plugging $\vy' = \vy_\gamma^*(\vx)$ into inequality (\ref{eq:grad-dom-last}) and using Cauchy-Schwarz inequality in the left-hand side gives 
  \begin{align}
        \Vert \nabla f_{\vx, \gamma}(\vy) \Vert \ge \frac{\gamma}{2} \Vert \vy - \vy_\gamma^*(\vx) \Vert^2.  \label{eq:cubic-grad-dom-1} 
    \end{align}
Moreover, inequality (\ref{eq:grad-dom-last}) also implies \citep[Lemma 4.2.1]{nesterov2018lectures} that
\begin{align*} 
    f_{\vx, \gamma}(\vy') \ge f_{\vx, \gamma}(\vy) + \langle \nabla f_{\vx, \gamma} (\vy), \vy' - \vy \rangle  + \frac{\gamma}{6} \Vert \vy' - \vy \Vert^3, \quad \forall \vy, \vy' \in \sR^d.
\end{align*}
Then, taking the minimum with respect to $\vy'$ on both sides of the above inequality gives
\begin{align} \label{eq:cubic-grad-dom-2}
f_{\vx, \gamma}^* \ge f_{\vx, \gamma}(\vy) - \frac{2}{3} \sqrt{\frac{2}{\gamma}} \Vert \nabla f_{\vx, \gamma}(\vy) \Vert^{\frac{3}{2}}.
\end{align}
Finally, inequality (\ref{eq:cubic-grad-uniform}) is merely the combination of (\ref{eq:cubic-grad-dom-1}) and (\ref{eq:cubic-grad-dom-2}).
\end{enumerate}
\end{proof}

\section{Proof of Lemma \ref{lem:sub-solover} }

\begin{proof}
According to Lemma \ref{lem:sub-prob-grad-dom} (item 1), the proximal function $f_{\vx,\gamma}$ has $(L+ 2 \gamma)$-Lipschitz Hessians. Based on this fact, we analyze the progress in the first, last, and other epochs separately. For brevity, we denote the suboptimality gap of each epoch as $F_{s}:= f_{\vx,\gamma}(\vy_s)  - f_{\vx,\gamma}^* $.
\begin{enumerate}
    \item At the beginning, Line \ref{lin:ini-CRN} in Algorithm \ref{alg:A-LEN} takes a CRN step on the function $f_{\vx,\gamma}$ from the initial point $\vx$. From \citep[(4.2.30)]{nesterov2018lectures} with $L_3 = M = L + 2 \gamma$, we have
    \begin{align} 
          f(\vy_0) \le \min_{\vy \in \sR^d} f_{\vx,\gamma}(\vy) + \frac{L+2\gamma}{3} \Vert\vy - \vx \Vert^3 
        \le f_{\vx,\gamma}^* + \frac{L+2\gamma}{3} \Vert\vy_\gamma^*(\vx) - \vx \Vert^3.
    \end{align}
    It means that
    \begin{align} \label{eq:gap-first-epoch}
        F_0 \le \frac{L+2\gamma}{3} \Vert\vy_\gamma^*(\vx) - \vx \Vert^3.
    \end{align}
    \item For the subsequent $S$ epochs ($s = 0,\cdots,S-1$), we show that choosing the parameters appropriately can reduce the suboptimality gap of each epoch by $1/2$. From inequality~(\ref{eq:final-bound-grad}), we have
    \begin{align} \label{eq:epoch-s-grad-small}
        \Vert \nabla f_{\vx,\gamma}(\vy_{s+1}) \Vert^{3/2} = \gO \left( \frac{m^{1/3} \sqrt{L+2\gamma} F_s}{N} \right).
    \end{align}
    In addition, from inequality (\ref{eq:cubic-grad-dom-2}), we also have 
    \begin{align} \label{epoch-s-grad-dom}
        F_{s+1} = \gO \left( \frac{\Vert \nabla f_{\vx,\gamma}(\vy_s) \Vert^{3/2}}{\sqrt{\gamma}
        } \right).
    \end{align}
    Hence, a combination of the above two inequalities gives 
    \begin{align}
        F_{s+1} = \gO \left( \frac{m^{1/3}  F_s}{N } \sqrt{\frac{L+2\gamma}{\gamma}} \right).
    \end{align}
    Therefore, if we set $ N = \Theta \left( m^{1/3} \sqrt{L/\gamma+2} \right) $, then we can ensure that $F_{s+1} \le F_s /2$ for all $s = 0,\cdots,S-1$, which further indicates that 
    \begin{align} \label{eq:gap-middle-epoch}
       F_S \le 2^{-S} F_0.
    \end{align}
    \item For the last epoch ($s = S$), we simply invoke inequality (\ref{eq:epoch-s-grad-small}) and the setting of $N$ to obtain
    \begin{align} \label{eq:gap-last-epoch}
    \Vert \nabla f_{\vx,\gamma}(\vy_{S+1}) \Vert^{3/2} = \gO \left( \sqrt{\gamma} F_S \right).
    \end{align}
\end{enumerate}
Finally, for $\hat \vy = \vy_{S+1}$, we combine the inequalities (\ref{eq:gap-first-epoch}), (\ref{eq:gap-middle-epoch}), and (\ref{eq:gap-last-epoch}) to obtain that
\begin{align*}
    \Vert \nabla f_{\vx,\gamma}(\hat \vy) \Vert^{3/2} = \gO \left( 2^{-S}
    \sqrt{\gamma} (L+2\gamma) \Vert \vy_\gamma^*(\vx) - \vx \Vert^3
    \right).
\end{align*}
If we choose $S = \Theta\left(\ln \left( \frac{1}{\sigma} \left( \frac{L+2\gamma}{\gamma} \right)^{2/3} \right)   \right )$, then we can ensure 
\begin{align*}
    \Vert \nabla f_{\vx,\gamma}(\hat \vy) \Vert \le \frac{\sigma \gamma}{6} \Vert \vy_\gamma^*(\vx) - \vx \Vert^2.
\end{align*}
From this inequality, we can further derive that
\begin{align*} 
       \Vert \nabla f_{\vx, \gamma}(\hat \vy) \Vert \le  \frac{\sigma \gamma}{3} \Vert \vx - \hat \vy \Vert^2 +  \frac{\sigma \gamma}{3}  \Vert \hat \vy - \vy_\gamma^*(\vx) \Vert^2 \le  \frac{\sigma \gamma}{3} \Vert \vx - \hat \vy \Vert^2 + \frac{2\sigma}{3} \Vert \nabla f_{\vx, \gamma}(\hat \vy) \Vert,
\end{align*}
    where the first inequality 
    uses the Young's inequality  $\Vert \va + \vb \Vert^2 \le 2 \Vert \va \Vert^2 + 2 \Vert \vb \Vert^2$, and the second one uses inequality
    (\ref{eq:cubic-grad-dom-1}). Then, rearranging the above inequality and using $\sigma \in (0,1)$ gives
    \begin{align} \label{eq:second-prox-condition}
     \Vert \nabla f_{\vx, \gamma}(\hat \vy) \Vert \le \frac{\sigma \gamma}{3 - 2 \sigma}  \Vert \vx - \hat \vy \Vert^2 \le \sigma \gamma  \Vert \vx - \hat \vy \Vert^2,
\end{align}
which is equivalent to the MS condition in Eq.~(\ref{eq:MS-oracle}).
\end{proof}

\section{Proof of Theorem \ref{thm:ALEN}}

\begin{proof}
Setting the parameter $\gamma = \Theta(L / m^{4/3})$ and picking any constant $\sigma \in (0, 0.99]$, we know that the parameters $S$ and $N$ becomes $S = \Theta( \ln m)$ and $N = \Theta(m)$. 
Moreover, invoking Lemma~\ref{thm:ALEN} with $\gamma = \Theta( L / m^{4/3})$, we know that Algorithm \ref{alg:A-LEN} finds an $\epsilon$-solution of Problem~(\ref{prob:main}) in $T = \mathcal{O}( m^{-8/21} \epsilon^{-2/7} )$ iterations. Therefore, and the total iteration complexity of Algorithm \ref{alg:A-LEN} is $TSN = \gO( (m + m^{13/21} \epsilon^{-2/7}) \ln m )$ as claimed by the theorem.
\end{proof}

\section{Proof of Lemma \ref{lem:sub-problem}}

\begin{proof}
We prove the result for the trust-region Newton subproblem (\ref{eq:tr-sub-prob}) and the CRN subproblem~(\ref{eq:CRN-prob}) individually.
\begin{enumerate}
    \item For the trust-region Newton subproblem (\ref{eq:tr-sub-prob}), we know from \citep[Corollary~7.2]{conn2000trust} that the pair $(\vh^*,\tau^*) \in \sB_D^d(\vzero) \times \sR_+$ satisfies both condition (\ref{eq:root-cond}) and the complementary slackness condition $ \tau^*( \Vert \vh^* \Vert - D) = 0$. Then, there are two cases: (1) $\vh^*$ is in the interior. (2) $\vh^*$ is on the boundary. For the first case, we have $\Vert \vh^* \Vert < D$,  $\tau^* = 0$ and thus $\vh^* = -\mA^{-1} \vb$. For the second case, we have  $ \Vert (\mA + \tau^* \mI_d)^{-1} \vb \Vert = D$, which means $\tau^*$ is the root of the function $\varphi^{\rm TR}(\tau) :=
    \Vert (\mA + \tau \mI_d)^{-1} \vb \Vert - D$.
    \item  For the CRN subproblem (\ref{eq:CRN-prob}), its first-order optimality condition is
    \begin{align} \label{eq:CRN-first-cond}
        \vb + \mA \vh^* + \frac{M \Vert \vh^* \Vert \vh^*}{2} = \vzero.
    \end{align}
    From \citep[Proposition 1]{nesterov2006cubic}, we know that it also satisfies the following second-order optimality condition
    \begin{align}  \label{eq:CRN-second-cond}
        \mA + \frac{M \Vert \vh^* \Vert}{2} \succeq \mO_d.
    \end{align}
    Then, by defining the variable $\tau^* =  M \Vert \vh^* \Vert /2$, we know that the pair $(\vh^*,\tau^*)$ satisfies condition (\ref{eq:root-cond}) and $\tau^*$ is the root of the function $\varphi^{\rm CRN}(\tau) = \Vert (\mA + \tau \mI_d)^{-1} \vb \Vert - 2 \tau/ M$.
\end{enumerate}
Note that the functions only differ in the second term. We further define their common nonlinear part in the first term as $\bar \varphi(\tau) = \Vert (\mA + \tau \mI_d)^{-1} \vb \Vert$ and show that the function $\bar \varphi(\tau)$ is convex and monotonically decreasing in $\tau$. We let $\mA = \mQ \mLambda \mQ^\top $ be the spectral decomposition of the matrix $\mA$, where $\mLambda = {\rm diag}(\lambda_1,\cdots,\lambda_d)$ and $\mQ \in \sR^{d \times d}$ is an orthogonal matrix. We further define $\tilde \vb = \mQ^\top \vb$ be the rotated vector in the eigenspace. Then we have 
\begin{align*}
    \bar \varphi(\tau) = \sqrt{ \sum_{i=1}^d \frac{\tilde \vg_i^2}{(\lambda_i + \tau)^2} },
\end{align*}
which is both convex and monotonically decreasing in $\tau$ when $\tau > \max\{ -\lambda_d, 0 \}$.
\end{proof}


\section{Results for Strongly Convex Problems} \label{apx:ALEN-SC}

In this section, we consider the setting when the objective function in Problem (\ref{prob:main}) is strongly convex.

\begin{asm} \label{asm:strongly-convex}
We assume that the function $f(\vx)$ is $\mu$-strongly convex, \textit{i.e.}, 
\begin{align*}
    f(\vx') \ge f(\vx) + \langle \nabla f(\vx), \vx' - \vx \rangle + \frac{\mu}{2} \Vert \vx'- \vx \Vert^2, \quad \forall \vx, \vx' \in \sR^d.
\end{align*}
\end{asm}

For strongly convex problems, we can target to find the following stronger solutions. 

\begin{dfn} \label{dfn:epsilon-root}
Under Assumption \ref{asm:strongly-convex}, we 
denote $\vx^* = \arg \min_{\vx \in \sR^d} f(\vx) $ as the unique solution to Problem (\ref{prob:main}). We call $\hat \vx$ an $\epsilon$-root if we have $\Vert \hat \vx - \vx^* \Vert^2 \le \epsilon$.     
\end{dfn}

Motivated by \citet{arjevani2019oracle,gasnikov2019optimal}, given a globally convergent algorithm for convex problems, we can apply a restart scheme to convert its convergence guarantee for finding an $\epsilon$-root in strongly convex problems in a black box manner. Following this recipe, we present the CALEN-restart method in Algorithm \ref{alg:A-LEN-SC}, which  repeatedly executes CALEN for a bounded number of steps, reducing the distance to $\vx^*$
by a constant factor in each epoch. We also formally state the main result of CALEN-restart for strongly convex problems as follows.

\begin{algorithm*}[htbp]  
\caption{CALEN-restart($f,\vx_0,S$)}  \label{alg:A-LEN-SC}
\begin{algorithmic}[1] 
\STATE \textbf{for} $s = 0,1,\cdots, S-1$ \\ [1mm]
\STATE \quad Select $T_s$ according to Eq. (\ref{eq:para-Ts}).
\STATE \quad $\vx_{s+1} = \text{CALEN}(f, \vx_s, T_s)$ \\[1mm]
\STATE \textbf{end for} \\[1mm]
\STATE \textbf{return} $x_{S}$  
\end{algorithmic}
\end{algorithm*}

    

\begin{thm} \label{thm:ALEN-SC}
Under Assumption \ref{asm:Hessian-Lip} and \ref{asm:strongly-convex}, setting the parameter $S = \lceil \log_2(R_0/\epsilon) \rceil $ and $\{T_s\}_{s=0}^{S-1}$ according to Eq. (\ref{eq:para-Ts}), then Algorithm \ref{alg:A-LEN-SC} can find an $\epsilon$-root in total $ \sum_{s=0}^{S-1} T_s = \gO( (m \ln (R_0/\epsilon) + m^{13/21} (L R_0/\mu)^{2/7}) \ln m )$ iterations, where $R_0 = \Vert \vx_0 - \vx^* \Vert$.

In particular, choosing $m = \Theta(\bar d)$ to minimize the total oracle costs in Eq. (\ref{eq:total-oracle-cost}), the equivalent gradient complexity of Algorithm \ref{alg:A-LEN} is $ \gO( (\bar d \ln (R_0/\epsilon) + \bar d^{13/21} (L R_0/\mu)^{2/7}) \ln \bar d)$.
\end{thm}

\begin{proof}
By the strong convexity of $f$, we know from \citep[Lemma 2.1.8]{nesterov2018lectures} that
\begin{align} \label{eq:SC-convert}
    \frac{\mu}{2} \Vert \vx - \vx^* \Vert^2 \le f(\vx) - f(\vx^*), \quad \forall \vx \in \sR^d.
\end{align}
We exploit this inequality to show that Algorithm \ref{alg:A-LEN-SC} with appropriate parameters 
can ensure that $\Vert \vx_s - \vx^* \Vert \le 2^{-s} R_0$ for all $s=  1,\cdots,S$. We prove it by induction. For brevity, we define $R_s : =\Vert \vx_s - \vx_0 \Vert $ and $F_s = f(\vx_s) - f(\vx^*)$. 
Suppose we have already ensured that $R_s \le 2^{-s} R_0$, we select parameters to ensure that $R_{s+1} \le 2^{-(s+1)} R_0$. From inequality (\ref{eq:SC-convert}), it suffices to let 
\begin{align} \label{eq:F-s1}
    F_{s+1} \le \frac{R_s^2}{8 \mu }.
\end{align}
Next, according to Theorem \ref{thm:ALEN} and $R_s \le 2^{-s} R_0$, the CALEN algorithm can ensure inequality~(\ref{eq:F-s1}) holds if we set
\begin{align} \label{eq:para-Ts}
\begin{split}
    T_s =& \gO\left( \left(m + m^{13/21} \left( L R_s/ \mu \right)^{2/7}  \right) \ln m  \right) \\
    =& \gO \left( \left( m  + m^{13/21} 2^{-2s/7} \left( LR_0/\mu \right)^{2/7} \right) \ln m \right).    
\end{split}
\end{align}
As this parameter setting ensures that $R_{s+1} \le 2^{-(s+1)} R_0$, we can finish the induction and obtain
\begin{align} \label{eq:restart-ALEN-s}
    R_S \le 2^{-S} R_0.
\end{align}
which means that the algorithm finds an $\epsilon$-root in $S = \lceil \log_2(R/\epsilon) \rceil$ epochs. Finally, the total iterations $\sum_{s=0}^{S-1} T_s$ can be bounded by Eq. (\ref{eq:para-Ts}) and the geometric series sum $\sum_{s=0}^{S-1} 2^{-2s/7} = \gO(1)$.

\end{proof}

\end{document}